\documentclass{amsart}
\usepackage{extarrows}
\usepackage{epsfig}
\usepackage{nnsy}
\usepackage{wasysym}

\def\ds{\displaystyle}
\def\R{\mathbb{R}}

\title{Limits and the System of Near-Numbers}
\author{Frank J.~Swenton\\ Middlebury College\\ Middlebury, VT}
\date{7 December 2004}
 
\begin{document}

\setlength{\floatsep}{10pt}
\setlength{\textfloatsep}{10pt}
\setlength{\intextsep}{10pt}
\setlength{\abovecaptionskip}{4pt}
\setlength{\belowcaptionskip}{10pt}

\begin{abstract}
In this paper, we present a comprehensive system for the treatment of
the topic of limits---conceptually, computationally, and formally.
The system addresses fundamental linguistic flaws in the standard presentation
of limits, which attempts to force limit discussion into the language of
individual real numbers and equality.  The system of \emph{near-numbers} properly fit the context 
of limits, allowing precise and unambiguous notation in limit computation.
More importantly, the near-numbers are a set of concrete mathematical objects
that provide a scaffold both for building an accurate intuitive understanding of
the limit concept and for connecting that understanding to the
formal definition of the limit.
\end{abstract}

\maketitle

\section{Introduction}
Within the first-semester calculus curriculum, the limit is the most pivotal and, simultaneously, the most troublesome concept for the student to grasp.  In this paper, we argue that a large number of the difficulties, both specific and general, that occur in the instruction of limits stem from the lack of a mathematical \emph{language} that properly addresses the fundamental nature of limits conceptually, computationally, and logically.  In particular, the notion of limit is
an operation that lacks a proper set of objects to operate on.  In response, we present an intuitive and logically sound mathematical language for use in the instruction of limits, which we call the system of \emph{near-numbers}.  

The near-numbers find their origin within elements of common limit notation and extend them into a robust, concise, and precise system within which to treat limits.  This system breaks the ambiguities in the standard notation, allows for a more gradual approach to the notions of limit and convergence, and provides additional structure both to the process of computing limits and to the standard logical definition of the limit---all while maintaining full formal correctness.

We begin below with a few statements in the language of near-numbers, the intents of which will be largely evident.  We will proceed to give a brief outline of the origin of the near-number system, initiating our discussion of 
the linguistic issues that it addresses and the pedagogical improvements 
it enables.  Before entering into a full discussion of the near-number system, we will give a quick trip through the near-numbers, including a sampling of specific pedagogical issues that it resolves within the curriculum.  Finally, we present the details of the near-number system and discuss some of the means by which
which it can enhance student comprehension of this troublesome topic in introductory calculus.  The appendix provides a brief presentation of the formal mathematical definitions of the objects and relations discussed herein.

\def\mbw{1.3in}
\def\mbwtwo{1.35in}

\vspace{8pt}
\noindent\makebox[\mbw]{\hfill\emph{Near-numbers:}\;\;}
\makebox[\mbwtwo]{$0^\NNp, 1^\NNm, 2^\NNpm$, $5^\NNpmd$, $4^\NNmd$\hfill}
%\makebox[\mbwtwo]{$2^\NNpd, 1^\NNmd, 5^\NNpmd$\hfill}
\makebox[\mbwtwo]{$\infty, -\infty, \pm\infty$\hfill}
$\NNstar$

\vspace{8pt}
\noindent\makebox[\mbw]{\hfill\emph{Arithmetic:}\;\;}
\makebox[\mbwtwo]{$1-\infty \to -\infty$\hfill}
\makebox[\mbwtwo]{$2^\NNp \times 3^\NNm \to 6^\NNpmd$\hfill}
$\ds\frac{-\infty}{\infty}\to\NNstar$

\vspace{8pt}
\noindent\makebox[\mbw]{\hfill\emph{Functions:}\;\;}
\makebox[\mbwtwo]{$\exp(-\infty) \to 0^\NNp$\hfill}
\makebox[\mbwtwo]{$\tan(\frac{\pi}2^\NNpm) \to \pm\infty$\hfill}
$\sin(\infty) \to \NNstar$

\vspace{8pt}
\noindent\makebox[\mbw]{\hfill\emph{Limits:}\;\;}
$\ds\lim_{x\to a} f(x) = L \Leftrightarrow f(a^\NNpm) \to L^\NNpmd$

\vspace{8pt}
\noindent\makebox[\mbw]{\hfill\emph{Continuity:}\;\;}
$f$ is continuous at $a \Leftrightarrow f(a^\NNpmd) \to f(a)^\NNpmd$

\vspace{8pt}
\noindent\makebox[\mbw]{\hfill\emph{Limit computation:}\;\;}
$\ds
\left.\frac1{\ln(\cos x)}\right|_{x\to\frac{\pi}2^\NNm}
\to \frac1{\ln(\cos(\frac{\pi}2^\NNm))}
\to \frac1{\ln(0^\NNp)} \to \frac1{-\infty} \to 0^\NNm$

%$\ds\lim_{x\to\frac{\pi}2^\NNm} \frac1{\ln(\cos x)}:\;\;\frac1{\ln(\cos(\frac{\pi}2^\NNm))}
% \to \frac1{\ln(0^\NNp)} \to \frac1{-\infty} \to 0^\NNm$

\section{The origin of the near-number system}

  	The origin of the near-number system is found in addressing linguistic flaws, both structural and semantic, of the standard notation for limits.  Whorf's Principle of Linguistic Relativity \cite{W} states, in short, that language impacts cognition at a fundamental level; we take this principle to heart in our approach:  the educational process is far more effective when both the student and the teacher are able to precisely express in writing the concepts at work in the topics at hand.  Near-numbers provide just such a language for limits; while the change they induce might at first appear merely cosmetic, we will see that a near-number presentation of limits enables deeper and more significant improvements in the overall pedagogy.

	The common expression $\lim_{x\to a}f(x) = L$ denotes an atomic logical assertion that the student is generally required to swallow whole, rather than having some well-defined, more basic mathematical ingredients from which to build it.  The two most fundamental flaws with this standard presentation of limits, and their consequent effect on the pedagogy, are apparent on re-examination.  The first is that the logical definition of $\lim_{x\to a}f(x) = L$  is too intricate to be instructed most effectively as a single atomic definition.  Second, the objects and notation that support this logical definition poorly fit the mathematics at work behind the scenes.  Equality between individual real
numbers is precisely what the limit concept is \emph{not} concerned with---thus both the objects (individual real numbers) and relation (equality) employed are in need of replacement.

	Moving away from the standard notation, it is common to speak more loosely: ``as $x\to a$, $f(x)\to L$'', or to be slightly more precise: ``as $x\to a^\NNpm$, $f(x)\to L^\NNpmd$''.  Here, the adornment on the $a$ indicates, as usual, that $x$ can be less than or greater than $a$, but not equal to $a$; in order that our notation be adequately precise, we add a dot to the adornment on $L$, indicating that $f(x)$ is permitted to equal $L$ (``on the dot,'' so to speak).  Our final notational observation is that the ubiquitous $x$ is entirely superfluous in this statement, which may be restated as ``$f(a^\NNpm) \to L^\NNpmd$''.  This more concise and precise notation for our original limit statement possesses a clear intuitive meaning of ``$f$ takes values near (but not equal to) $a$ to values near $L$.''  While this notation does indeed fix some of the flaws of imprecision inherent in the original statement, we have thus far merely cleaned it up and clarified it, but we have not yet made any significant structural changes.

Our next concern is that, thus far, the statement $f(a^\NNpm) \to L^\NNpmd$ is still atomic, rather than being (as it should appear) a compound statement built from simpler parts.  Clearly, the first necessary step is to treat $a^\NNpm$ and $L^\NNpmd$ as mathematical objects that exist in their own right, rather than being merely symbols tied to the limit notation---such objects are what we call \emph{near-numbers}, which form the cornerstone of our new language for use in the discussion of limits.  Next, we need real-valued functions (as well as arithmetical operations) to be capable of operating on near-numbers.  Third, we must complete this language by defining the relation ``$\to$'' between near-numbers.  Throughout, we must ensure that these concepts can be fully formalized, so that we are not merely working within a baseless intuitive system.  Via near-numbers, we will be able to parse the statement $f(a^\NNpm) \to L^\NNpmd$ as: ``start with the near number $a^\NNpm$, then apply $f$ to it; the result fits the near-number $L^\NNpmd$''.  As such, we will come to an understanding of its meaning---both intuitively and formally---one component at a time.

The system proposed in this paper does not insist upon any fundamental change to the logical definition of the limit.  It merely fleshes out those notions that anyone familiar with limits has conceptualized, providing the student with a scaffold from which to build intuitive understanding and computational skills (as well as the formal logical definitions that support them, to whatever extent this is desired).  
The path provided by near-numbers will result in the standard logical definitions of limit statements---in fact, the near-number approach will lead us directly to a seamless treatment of all finite and infinite limit definitions simultaneously.

\section{A quick trip through the near-numbers}

Before entering into a full discussion of near-numbers, we first give a quick overview of the near-number system and a few of its pedagogical features.  It should be noted that this notation is fully rigorous, as each statement made herein with near-numbers can be formally defined and proven.  The near-numbers themselves are largely recognizable in their intent, for they find their roots in standard limit notation; the primary difference is that these symbols retain intrinsic meaning even outside the limit context.  The basic near numbers fall into three categories: the \emph{finite} near-numbers ($0^\NNp$, $1^\NNpm$, $L^\NNpmd$, etc.---a dot indicates that the number itself is to be considered, which is otherwise not the case); the \emph{infinite} near-numbers ($+\infty$, $-\infty$, $\pm\infty$); and the \emph{indeterminate} near-number ($\NNstar$).  

Near-numbers can be acted upon by the standard arithmetic operators and real-valued functions, but the equality symbol is discarded.  Instead, we have the relation ``$\to$'' between near-numbers, which expresses the notion of one near-number ``fitting into'' another; unlike equality, this relation has a crucial property of directionality (it is reflexive and transitive, but not symmetric).  The assertion $\alpha\to\beta$ is logically defined as true or false for any pair $\alpha,\beta$ of near-numbers; for example, it is true that $2^\NNp + 3^\NNm \to 5^\NNpmd$, $\cos(0^\NNpm) \to 1^\NNm$, and $2^\NNp \to 2^\NNpm$, but $2^\NNpm \to 2^\NNp$ is false.  

Expressed in terms of near-numbers, basic limit statements become much cleaner and more concise.  For example, $\lim_{x\to a} f(x) = L$ becomes simply $f(a^\NNpm) \to L^\NNpmd$ [``$f$ maps values near (but different from) $a$ to values near $L$''], and the definition of a real-valued function $f$ being continuous at $a$ changes from $\lim_{x\to a}f(x) = f(a)$ to the more fluid $f(a^\NNpmd) \to f(a)^\NNpmd$ [``$f$ maps values near $a$ to values near $f(a)$''].
Below, we present three simple examples of the shortcomings of the standard limit language, discussing each within the context of near-numbers to illustrate just a few of the logical and pedagogical advantages of the system.

\vspace{8pt}
\noindent\emph{Example 1.\;\;}
$\ds\frac10$ versus $\ds\lim_{x\to0}\frac1{x^2}$ versus $\ds\lim_{x\to0}\frac1{x^3}$
\vspace{4pt}

In the case of $\frac10$, division by zero is undefined.  Yet we would say that as $x\to0$, $x^2\to0$, so that attempting to work structurally within the real number system, one operation at a time, we find that the two limits simplify to $\frac10$ as well.  Certainly, we would argue that in the latter cases, the zero \emph{does not mean zero}---but given this fact, we should demand a more precise language that allows us express what we \emph{do} mean.  Such a shortcoming of our limit language forces us into an unnecessary symbolic ambiguity, resulting in confusion and even apparent contradiction in the mind of the student.  In light of Whorf's Principle of Linguistic Relativity, this is an issue in serious need of address.
	
In the near-number system, the first quantity $\frac10$ is still undefined (the laws of real numbers remain intact).  However, we can express the limit computations more precisely via near-numbers, substituting for $x$ and simplifying one operation at a time.  We can write these computations, showing each step, as
follows:

\vspace{6pt}
\hfill
$\ds \left.{\frac1{x^2}}\right|_{x\to0^\NNpm}\to \frac1{(0^\NNpm)^2} \to \frac1{0^\NNp} \to +\infty$;
\hspace{0.2in}
$\ds \left.{\frac1{x^3}}\right|_{x\to0^\NNpm}\to \frac1{(0^\NNpm)^3} \to \frac1{0^\NNpm} \to \pm\infty$
\hfill
\vspace{6pt}

Here, what is written more clearly expresses what is meant (indeed, the exact real number zero occurs logically not once in these two limits---so why should it appear symbolically?).  The limits can be computed structurally, and these three situations that typically present an apparent contradiction due to symbolic ambiguity are distinguished in the intermediate work as $\frac1{0}$, $\frac1{0^\NNp}$, and $\frac1{0^\NNpm}$.

\vspace{8pt}
\noindent\emph{Example 2.\;\;}
If $\ds\lim_{x\to a}g(x) = b$ and $\ds\lim_{y\to b}f(y) = L$, 
then $\ds\lim_{x\to a} f(g(x)) = L$.  \quad[False!]
\vspace{4pt}

The shortcomings of the standard limit language are very clearly demonstrated in this well-known example of a statement that appears true to the beginner (and often, at first glance, to one who understands limits quite well).  Not only does the notation obscure the issue at hand due to its imprecision, but it is symbolically clumsy with its wash of $\lim$'s and dummy variables.

In contrast, the above statement simplifies greatly and becomes discernibly false when expressed via near-numbers.  Recalling the formulation of continuity via near-numbers, the two hypotheses become $g(a^\NNpm)\to b^\NNpmd$ and $f(b^\NNpm)\to L^\NNpmd$, and the conclusion becomes $f(g(a^\NNpm))\to L^\NNpmd$.  Here, we can attempt to simplify the left-hand side of the conclusion inside-out, as usual.  Our first hypothesis tells us that $f(g(a^\NNpm)) \to f(b^\NNpmd)$; from this point, it not only becomes clear that we can make no further conclusion (our second hypothesis is on $f(b^\NNpm)$, not $f(b^\NNpmd)$!), but it also makes clear the precise cause of the logical problem with this assertion---namely, that we do not have any information about $f(b)$ itself.

\vspace{8pt}
\noindent\emph{Example 3.\;\;}
If $f$ is continuous at $a$ and $g$ is continuous at $f(a)$, 

\noindent\phantom{\emph{Example 3.\;\;}}
then $g\circ f$ is continuous at $a$.
\vspace{4pt}

As is evident from consideration of the previous remark, it is not possible to write a formally correct justification of the above assertion using only standard limit notation and the definition that $f$ is continuous at $a$ if $\lim_{x\to a}f(x)=f(a)$  (one, to be correct, is forced to resort to $\varepsilon$-$\delta$ definitions).  However, via the near-numbers our definition of continuity of $f$ at $a$ is neatly phrased as $f(a^\NNpmd)\to f(a)^\NNpmd$.  Thus, the proof of the assertion becomes a trivial matter of pulling out $\NNpmd$'s: $g(f(a^\NNpmd)) \to g(f(a)^\NNpmd) \to g(f(a))^\NNpmd$, where the two steps are justified by the first and second hypotheses, respectively.

\section{The Near-Number System}
\subsection{Intent}

Our proper introduction to near-numbers begins at the intuitive level, establishing what it is that we are \emph{trying} to express with a near-number.  Once this foundation has been laid, we'll make these notions more precise, leading eventually to the mathematical definitions of these objects.  To begin, we pose and answer the two key questions about near-numbers:

\def\mbw{1.7in}

\vspace{6pt}
\hangindent=1.7in
\noindent\makebox[\mbw]{\hfill\emph{What is a near-number?}\;\;}%
A near-number is \emph{not} a number (nor even a set of numbers), but an object representing the numbers \emph{near} some location on the real line.  

\vspace{4pt}

\noindent\makebox[\mbw]{\hfill\emph{How near?}\;\;}%
Nearer than any particular real number.
\vspace{6pt}

\noindent The utility of near-numbers derives from the simple fact that when considering limits of any form, it is never individual numbers that are being considered.  The near-number system provides a set of objects via
which we can properly discuss such concepts.

\begin{figure}[!h]
\begin{centering}
\epsfbox{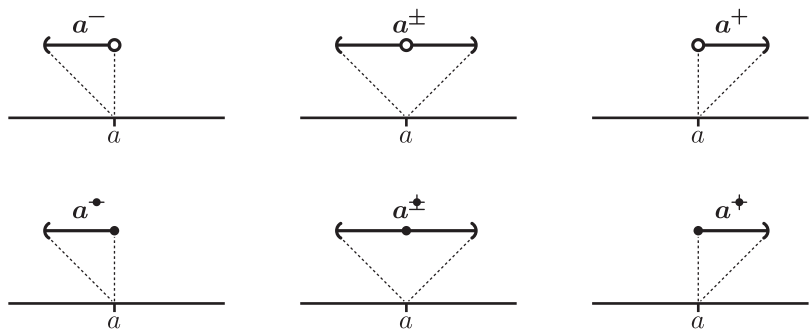}
\end{centering}
\caption{The finite near-numbers}
\label{Finite}
\end{figure}

The basic near-numbers fall into three categories:  the \emph{finite near-numbers}, the \emph{infinite near-numbers}, and the \emph{indeterminate near number}.  The finite near numbers are written by adorning a real number with some symbol indicating which nearby numbers are under consideration.  The first three of these adornments are already familiar: $a^\NNp$, $a^\NNm$, and $a^\NNpm$ represent, respectively, the real numbers ``just to the right of'', ``just to the left of'', and ``just to the left and right of'' the real number $a$ (as usual, the value $a$ itself is not included in these near-numbers).  Because situations commonly arise in which we need to allow the number $a$ itself (a subtle but important distinction), we provide for this case: the addition of a small dot or circle to the adornment of any finite near-number (e.g., $a^\NNpd$, $a^\NNmd$, $a^\NNpmd$) indicates that the number $a$ itself is to be included.  We can picture these six types of finite near-numbers as the ``infinitesimal intervals'' in the real line, as indicated in Figure~\ref{Finite}.

These near-numbers encapsulate and give a body to the precise notions with which we attempt to work every time that we consider finite limits, and we will see that they give a logically valid meaning to notions such as ``the numbers just after $a$.''  While this might be considered to be dangerous, the two fundamental principles of near-numbers avert such danger.  First, a near-number is not a number; fully impressing upon the student that they are starting from scratch with near-numbers and that, \emph{a priori}, none of the rules they've come to expect from numbers (cancellation, arithmetic, equality, etc.)~apply, then there is no existing context into which near-numbers should intrude to cause any problems.

Second, we might note the apparent lack of mathematical precision in saying ``near''; but again, the principle ``nearer than any particular real number'' clearly impresses what is demanded.  Does ``$0^\NNp$'' mean ``$0.00001$''?  No, for it represents numbers nearer to zero than any particular real number---in particular, nearer than $0.00001$ (indeed, we can picture the point in the diagram at which the near-number excludes any given positive real number).  Does this cause a problem?  No, for $0^\NNp$ is not a real number, nor is it even a particular set of real numbers.  

\begin{figure}[!h]
\begin{centering}
\epsfbox{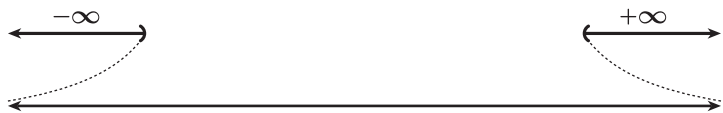}
\end{centering}
\caption{The infinite near-numbers $-\infty$ and $+\infty$; $\pm\infty$ is their union}
\label{Infinite}
\end{figure}

The three infinite near-numbers are $+\infty$ (also denoted by $\infty$), $-\infty$, and $\pm\infty$;  these indicate, respectively, the numbers ``near the right end,'' ``near the left end,'' and ``near either end'' of the real line.  As usual, the answer to ``how near?'' is ``nearer than any particular real number''---this treatment of proximity will allow us to smoothly handle limits and limit definitions for both finite and infinite cases in a uniform manner, leading to a better conceptualization of the notions.  As with the finite near-numbers, we can represent the intent of these three near-numbers via diagrams, as in Figure~\ref{Infinite}.

The primary difference between the finite near-numbers and these infinite near-numbers is that proximity is measured not by distance, but by order:  to be ``nearer the right end'' than $M$ means to be greater than $M$.  If the real line were to have left and right endpoints, this notion would directly parallel that of proximity to a finite near-number, so these infinite near-numbers have the same basic character as the finite near-numbers.  As with the finite near numbers, we can pose the question: ``does $+\infty$ mean $10000$?''  The answer is that, no, it represents larger numbers---numbers larger than any particular real number, as again is visibly indicated on the diagram.

The one remaining near-number is the \emph{indeterminate near-number}, $\NNstar$, which serves as a catch-all for any situation that is unknown, is not specified by one of the other near-numbers, or cannot be logically determined.  We declare computation involving a $\NNstar$ to be illegal;  in this sense, the near-number $\NNstar$ indicates failure of a computation, i.e., it constitutes a stop sign, or something of a computational black hole.  Given that this special near-number indicates nothing whatsoever about a position on the real line, we choose not to diagram it (though this indeterminacy could be represented by taking $\NNstar$ to fill the entire real line).

\subsection{Near-number diagrams revisited}

\begin{figure}[!h]
\begin{centering}
\epsfbox{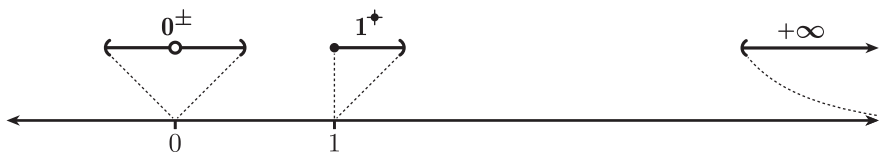}
\end{centering}
\caption{Infinitesimal view of $0^\NNpm$, $1^\NNpd$, and $+\infty$}
\label{Infinitesimal}
\end{figure}

While the idea of near-numbers being ``infinitesimal intervals'' (as Figure~\ref{Finite} seems to indicate) can be quite intuitively seductive, we should require a more concrete logical formalization of these objects before we accept them as fit for mathematical use.  In this hope, we are fortunate that the above intuitive diagrams provide just what we seek when read properly.  Our motivation is as follows: when we shrink our sets as we follow the dotted lines downward, nothing significant remains at the end of the process---so our interest must properly be in the progression itself and not simply the ``bottom line''.
We illustrate our discussion via the examples of the three near-number diagrams in Figure~\ref{Infinitesimal}.

\begin{figure}[!h]
\begin{centering}
\epsfbox{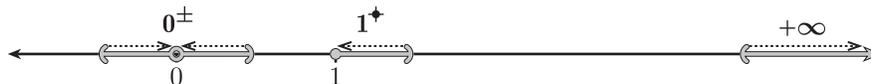}
\end{centering}
\caption{Animated-set view of $0^\NNpm$, $1^\NNpd$, and $+\infty$}
\label{Animated}
\end{figure}

Every near-number diagram can be profitably interpreted in two alternative manners, each of which expresses the near-number in a slightly different light.  One perspective is obtained by considering each near-number diagram to give us instructions for an animation of sets in the real line, starting from the top of the infinitesimal-view diagram and working toward the bottom, as in Figure~\ref{Animated}.  

\noindent%
Each frame in such an animation expresses the sense of being progressively
``nearer'' to the concept expressed by the near-number.  As such, we can think
of the entire animation to fully express what is demanded in order to have
the properties of, say, ``$\infty$'' or ``$1^\NNpd$'': to become arbitrarily close to some
location on the real line---even if the exact location is unattainable.
With modern computer technology being easily capable of displaying such animations, this presentation is very accessible.

\begin{figure}[!h]
\begin{centering}
\epsfbox{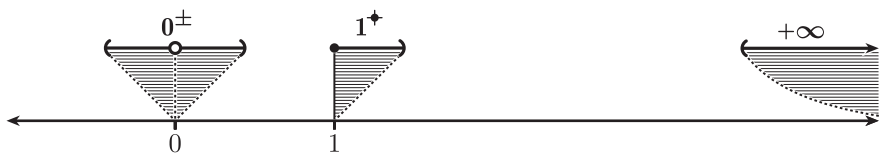}
\end{centering}
\caption{Static planar-set view of $0^\NNpm$, $1^\NNpd$, and $+\infty$}
\label{Static}
\end{figure}

A third perspective on near-numbers is obtained by dispensing with time and animation, instead considering each near-number diagram to express a set in the plane, as in Figure~\ref{Static}.  It is useful to view these sets as collections of horizontal ``slices,'' which keeps with the spirit of the animated view of near-numbers (the slices correspond to the frames of the animation).  
Again, each slice paints for us a target expressing what
it means to be close to the near-number concept at hand; the primary difference
is that in this view we see all of the targets at once.
Because this perspective puts the entirety of a near-number in view simultaneously, it will be the most useful perspective to take when discussing mathematical definitions and formal proof.  

Thus, we can think about near-numbers intuitively as ``infinitesimal intervals'', dynamically as an animated progression of sets, or geometrically via static sets in the plane.  Each of these perspectives has its own benefits; while we often intuitively think of near-numbers via the first view, we will generally prefer the second or third when we wish to be more precise.  These three interrelated perspectives provide a great deal of pedagogical flexibility, allowing the presentation of near-numbers to be molded to the needs of a variety of contexts.

\subsection{The arrow relation}
\label{Arrow}

Within the near-number system, the symbol ``$\to$'' represents a certain relation among near-numbers; unlike the relation of equality for the real numbers, the arrow is necessarily directional---reflexivity and transitivity hold, but symmetry generally does not.  The statement ``$\alpha\to\beta$'' can be read intuitively as ``$\alpha$ fits the form of $\beta$'' or ``$\alpha$ fits into $\beta$.''  We note that the common term ``goes to'' is entirely inappropriate, as nothing is going anywhere---there are two near-number $\alpha$ and $\beta$, and $\alpha$ either does or does not fit into $\beta$.

We can consider the arrow relation loosely on general principles for the basic near-numbers, simply interpreting their meanings or inspecting the ``infinitesimal'' near-number diagrams to obtain the relationships in Figure~\ref{ArrowRels} (note that none of the arrows are reversible).

\begin{figure}[!h]

\def\gapspace{0.2in}
\def\endspace{0.2in}

\def\gapbox#1{\makebox[\gapspace]{\hfill #1 \hfill}}
\def\endbox#1{\makebox[\endspace]{\hfill #1 \hfill}}

\hbox{
\hspace{0.4in}
\endbox{$a^\NNm$}\gapbox{$\to$}\endbox{$a^\NNpm$}\gapbox{$\leftarrow$}\endbox{$a^\NNp$}\hfill
}

\hbox{
\hspace{0.4in}
\endbox{$\downarrow\;\,$}\gapbox{}\endbox{$\downarrow\;\,$}\gapbox{}\endbox{$\downarrow\;\,$}
\hspace{0.5in}
\def\endspace{0.25in}
\endbox{$-\infty$}\gapbox{$\to$}\endbox{$\pm\infty$}\gapbox{$\leftarrow$}\endbox{$+\infty$}
\hspace{0.5in}
$\alpha \to \NNstar$
\hfill
}

\def\endspace{0.2in}

\hbox{
\hspace{0.4in}
\vspace{2pt}
\endbox{$a^\NNmd$}\gapbox{$\to$}\endbox{$a^\NNpmd$}\gapbox{$\leftarrow$}\endbox{$a^\NNpd$}
\vspace{6pt}
\hfill
}
\caption{The arrow relation for the basic near-numbers}
\label{ArrowRels}
\end{figure}

\noindent%
While our intuitive understanding suffices for these basic near-numbers, we will need a stronger formulation of the arrow, particularly once we begin to operate on near-numbers via arithmetic and real-valued functions.  Thus, we now
shift to our alternative perspectives on near-numbers.

In terms of the animated-set perspective on near-numbers, we can interpret the statement $\alpha\to\beta$  as ``$\alpha$ fits the form of $\beta$''.  First viewing $\beta$, we see a progression of sets in the real line.  This animation expresses the qualitative meaning of $\beta$, in terms of what a progression of sets must do to fit the its form; e.g., we could describe $0^\NNpm$ as missing zero but pinching in arbitrarily close to zero as the slices progress.  Now, viewing the animation of the near-number $\alpha$, we can answer whether we see the required qualities.

While such a qualitative view of the arrow relation expresses quite well what we intend by the arrow relation intuitively, we can state the definition of $\alpha\to\beta$ more specifically as follows:  each slice of $\beta$ paints a progressively smaller target on the real line; for each such target, there must be some slice of $\alpha$ that fits within that target.  This view directly reflects the standard sense of limit, and it avoids the ambiguities in the qualitative animated description of the arrow relation given above.  In particular, it is insignificant whether \emph{every} slice of $\alpha$ fits into some slice of $\beta$, and the rate at which $\alpha$ shrinks is irrelevant---all that matters is whether it eventually fits into every target painted by $\beta$.

We can readily translate the above quantitative description of the relation $\alpha\to\beta$ into the static planar-set view of near-numbers.  In this perspective, we will view all frames of the animation simultaneously, affording a more complete view of the situation and a more solid foundation from which to write formal definitions.  Considering each near-number's shape to be formed from horizontal slices, we can consider the process of ``dropping'' the slices of $\alpha$ vertically into $\beta$ and checking whether they can reach to the bottom of the shape;  if so, then $\alpha\to\beta$ (whence ``$\alpha$ fits into $\beta$''), and if not, then $\alpha\not\to\beta$.  

\begin{figure}[!h]
\begin{centering}
\epsfbox{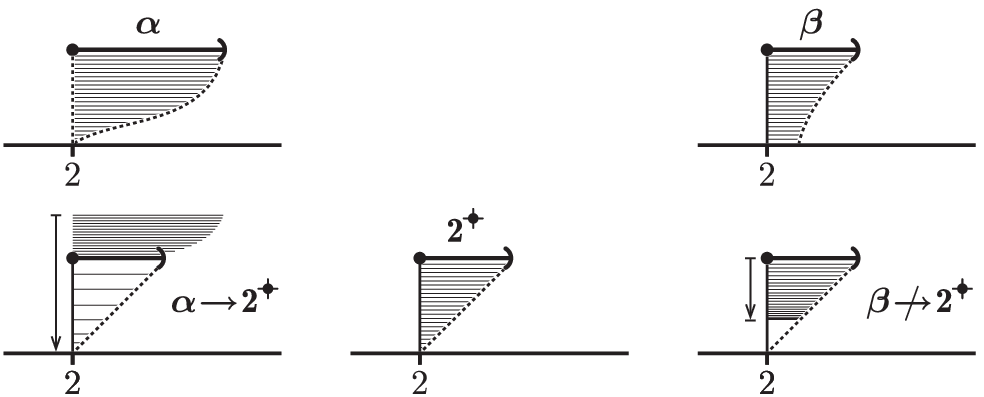}
\end{centering}
\caption{The arrow relation:  $\alpha\to2^\NNpd$, while $\beta\not\to2^\NNpd$}
\label{ArrowFig}
\end{figure}

Lining up two near-number diagrams, this method gives a simple geometric description of the arrow relation, as illustrated in Figure~\ref{ArrowFig}.

\noindent%
It is important to keep in mind
that this geometric construction is merely a simple and intuitive means
for expressing the same concept with which we have been concerned all along in 
the context of near-numbers.
As we move downward through the slices of $\beta$, we are presented an
increasingly refined view of what $\beta$ expresses about position
within the real line.  For the slices
of $\alpha$ to reach to the bottom of $\beta$ simply means that $\alpha$ fits
this description.

This geometric definition of the arrow relation, when quantified, will yield the standard formal definition of the limit for all finite and infinite cases, as will be discussed in Section~\ref{Quantify}.  Moreover, the details of the definition will emerge naturally from our simple geometric description, so that the order of quantification will not be in need of rote memorization. Thus, we obtain a natural path to what is generally the most problematic logical definition within the calculus curriculum.  For now, our geometric definition of the arrow relation will suffice, so we will discuss functions, arithmetic, and ``limits'' before moving on to the formal logical definition of the arrow relation.

\subsection{Real-valued functions}

\begin{figure}[!h]
\begin{centering}
\epsfbox{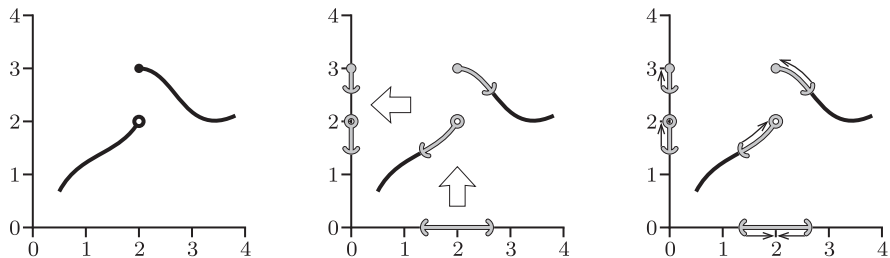}
\end{centering}
\caption{Mapping a set and a near-number (animated view)}
\label{Function}
\end{figure}

Given a near-number $\alpha$ and a real-valued function $f$, we can easily obtain the near-number $f(\alpha)$ via the animated view of the near-number $\alpha$, as follows.  Graphing $f$ and plotting a slice of $\alpha$ on the $x$-axis, we first use the graph of $f$ to map this set, just as we use it to map individual values:  starting with a set on the $x$-axis, we locate the corresponding set of points on the graph, then project this set to the $y$-axis, as in Figure~\ref{Function}.
If any value in a slice does not lie within the domain of $f$, we indicate the invalid operation via a $\NNstar$.  Considering the sets that result as we progress through the slices of the near-number $\alpha$, we obtain a precise animation of the near-number $f(\alpha)$.  Once we have obtained $f(\alpha)$, we are immediately able to address the question of whether $f(\alpha)\to\beta$ for a given near-number $\beta$, by recognizing the shape of the near-number result that we see.  Continuing with the example illustrated above, we can shrink our domain interval to find that $f(2^\NNpmd)\to 2^\NNm \cup 3^\NNmd$ or, more conventionally, that $f(2^\NNm)\to 2^\NNm$ and $f(2^\NNpd)\to 3^\NNmd$.  The rightmost two diagrams above, while fairly straightforward, are better envisioned as
the animations that they are meant to indicate;  when viewed as such, the
arrows fall away and the operation becomes quite trivial for any student to understand.

This setwise treatment of functions handles the limit operation for functions more properly than can any statement about $x$ ``going to'' $a$.  The near-number approach respects the fact that the limit has nothing to do with any particular value of $x$, allowing a clearer view of what the limit is concerned with and opening the full process to inspection.  In addition, near-numbers possess a \emph{monotonicity} that simplifies limit considerations.  For example, the oscillation of the functions $f(x)=\sin x$ and $g(x) = \frac1x\sin x$ ``as $x\to\infty$'' runs a serious risk of confusion in a student's mind over whether the second function indeed ``gets closer and closer to zero,'' given that when we consider one value at a time, the approach is not monotonic.  However, this issue (which is irrelevant to the limit) is not in play when we consider $f(\infty)$ and $g(\infty)$ via near-numbers, for the setwise action is monotonic and there is no room for confusion that $f(\infty)\to\NNstar$ (or, more precisely, $f(\infty)\to[-1,1]$), while $g(\infty)\to0^\NNpmd$.

\begin{figure}[!h]
\begin{centering}
\epsfbox{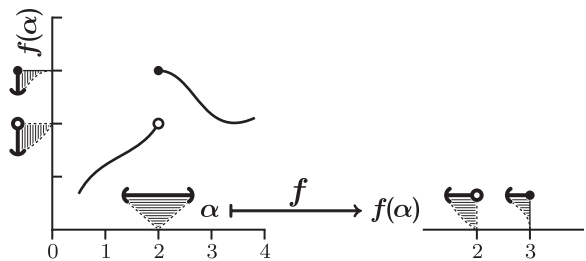}
\end{centering}
\caption{Mapping a near-number (static-set view)}
\label{Stack}
\end{figure}

We can directly translate our animated view of mapping near-numbers into the static planar-set perspective by (as usual) stacking our slices.  
Again, this operation is best viewed via an animation, which these pages are malequipped to provide.  Figure~\ref{Stack} indicates the result
of this stacking of slices, continuing from the rightmost diagram in
Figure~\ref{Function} to record what we observe in the animation.
Due to the fact that the graph of a function $f:\R\to\R$ poorly respects the cohabitation of the domain and range of $f$, our slices end up on the
vertical axis.  In order to most efficiently move such a set into its ``usual'' position on the real line, we stack the slices to the \emph{left side} of the $y$-axis so that we can put them into position via clockwise rotation about the origin.

Via this procedure, we can compute $f(\alpha)$ as a near-number in its own right, 
so that the statement $f(\alpha)\to\beta$ can be treated just like any other case of the $\to$ relation.  We note in closing that the values 
resulting for each slice of $f(\alpha)$ are simply the values $f(x)$ where $x$ is in the corresponding slice of $\alpha$, despite our chase through the Cartesian plane in obtaining them---we have only taken the graph-based approach in 
keeping with its ubiquity within the curriculum.  Properly using the notion of real-valued 
function, we can (and should) construct $f(\alpha)$ directly by applying
$f$ to the values in each slice of $\alpha$.

\subsection{Arithmetic}

The arithmetic operations of addition, subtraction, multiplication, and division can, as usual, be viewed from the standpoint of any of our three perspectives on near-numbers.  As a preliminary topic of discussion, we have arithmetic between a real 
number and a near-number.  Fortunately, our geometric interpretations of
real arithmetical operations remain valid:  adding or subtracting a real number
from a near-number shifts the near-number's position on the real line; and
multiplying or dividing a near-number by a real number scales the near-number
relative to the origin.  

We can summarize the arithmetic between near-numbers via rules for addition, negation, multiplication, and reciprocation in just a few lines each:

\def\mbw{1in}
\def\myrightbox#1{\noindent\hangindent=\mbw \makebox[\mbw]{\hfill #1\;\;}}
\def\myleftbox#1{\noindent\hangindent=\mbw \makebox[\mbw]{#1 \hfill}}

\vspace{8pt}
\myleftbox{ADDITION}

\vspace{2pt}
\myrightbox{\emph{Special:}}%
$-\infty + \infty \to \NNstar$; otherwise, $\infty+\alpha\to\infty$ and $-\infty+\alpha\to-\infty$.

\vspace{2pt}
\myrightbox{\emph{General:}}%
Add the base numbers.  For the adornment, any $\NNp/\NNm$ conflict in the summands' adornments gives $\NNpmd$; otherwise, the adornment is $\NNp$ or $\NNm$, with a dot just when \emph{both} have dots.

\vspace{8pt}
\noindent\myleftbox{\emph{NEGATION}}%
Negate both the base and the adornment.

\vspace{8pt}
\noindent\emph{MULTIPLICATION}\quad
[Factor out any negatives first; split $0^\NNpm$ and $\pm\infty$]

\vspace{2pt}
\myrightbox{\emph{Special:}}%
$\infty\times0^{?}\to\NNstar$; otherwise $\infty\times\alpha\to\infty$.

\vspace{2pt}
\myrightbox{\emph{General:}}%
Multiply the base numbers.  Adornment rules as for addition, except that a dotted zero in a factor always yields a dotted zero in the result.

\vspace{8pt}
\noindent\emph{RECIPROCATION}

\vspace{2pt}
\myrightbox{\emph{Special:}}%
Pairs of reciprocals: $+\infty$ and $0^\NNp$; $-\infty$ and $0^\NNm$; $\pm\infty$ and $0^\NNpm$.

\myrightbox{}%
$\ds\frac1{0^\NNpd},\frac1{0^\NNmd},\frac1{0^\NNpmd}\to\NNstar$ (division by zero!).

\vspace{4pt}
\myrightbox{\emph{General:}}%
Reciprocate the base and negate the adornment.

\vspace{6pt}

\noindent%
These rules can be memorized or reasoned out from the intuitive definitions, and they yield the rules of subtraction and division (which, being by nature composite, are a bit hairier).  We need not rely on such reasoning, however, because each of these results can be accurately determined visually from near-number diagrams.

In both the animated and static perspectives on near-number diagrams, we can view the above arithmetical operations geometrically to logically determine their results for any near-number arguments.  In either perspective, we view the near-numbers one slice at a time, i.e., via sets of real numbers.  We then operate on these sets of real numbers as usual and observe the resulting sets.  While such set operations are typically not a part of the first-year calculus curriculum, the topic can be treated in a very straightforward manner for the specific cases with which we'll be concerned.  Such a setwise approach to these operations is the only proper path to a clear view of indeterminacy, as we will see below.

\begin{figure}[!h]
\begin{centering}
\epsfbox{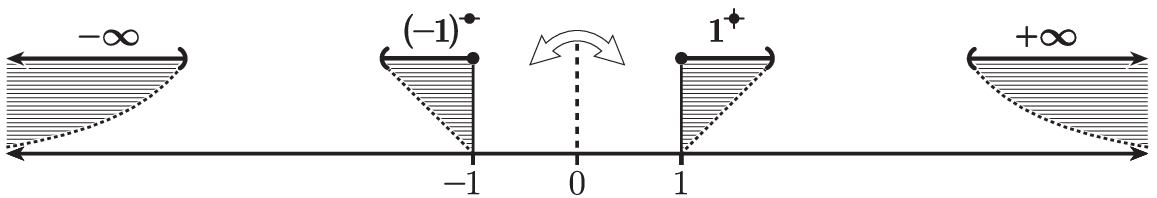}
\end{centering}
\caption{Negation:  
$-\infty\xleftrightarrow{\mathrm{negate}}+\infty$ 
and 
$(-1)^\NNmd\xleftrightarrow{\mathrm{negate}}1^\NNpd$}
\label{Neg}
\end{figure}

\begin{figure}[!h]
\begin{centering}
\epsfbox{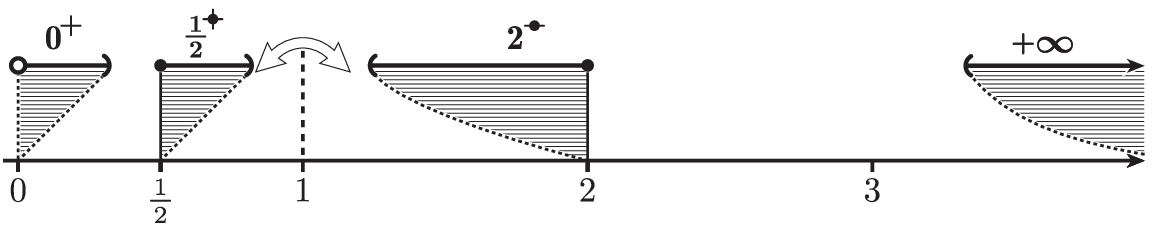}
\end{centering}
\caption{Inversion:  
$0^\NNp\xleftrightarrow{\mathrm{invert}}+\infty$
and 
$\frac12^\NNpd\xleftrightarrow{\mathrm{invert}}2^\NNmd$}
\label{Inv}
\end{figure}

In the cases of the unary operations of negation and reciprocation, we can very easily work with near-number arithmetic geometrically, interpreting these operations as mirroring the real line about zero and inversion of the real line about zero, respectively.
When the number zero itself is inverted, independent of any other concerns, the result is $\NNstar$ (this is only a concern when inverting the three near-numbers $0^\NNpd$, $0^\NNmd$, and $0^\NNpmd$).  These operations are simple to perform in both the animated view and the static view 
of near-numbers, so we present the more complete static view in Figures~\ref{Neg} and \ref{Inv}. Note that, as we are not concerned with equality (only arrow-equivalence), all shapes appearing are labeled with arrow-equivalent basic near-numbers.

When working with the binary operations of addition and multiplication (and subtraction and division), we consider a pair of slices, one from each parameter; without loss, these can be
taken at the same level.  With slices in hand, we must determine the set of all resulting sums or products.  This definition---that we're considering all possible sums or products of elements in one slice with those in the other---is fundamental and should not be overshadowed by any methods for determining it.  One method for this procedure is to, for
each element of one slice, perform the arithmetic operation on the other slice
(which amounts merely to translation or scaling).  Collecting
all such results as the chosen element ranges through the first slice yields a slice of the resulting near-number, which we can view either by animating 
these slices or by stacking them to obtain a static set.  Either view is best illustrated
by an animation showing the action performed, which again these pages are malequipped to provide.  We simply present a few illustrative
examples, paying particular attention to indeterminate cases, whose indeterminacy is very clearly exhibited via the near-number approach.

\begin{figure}[!h]
\begin{centering}
\epsfbox{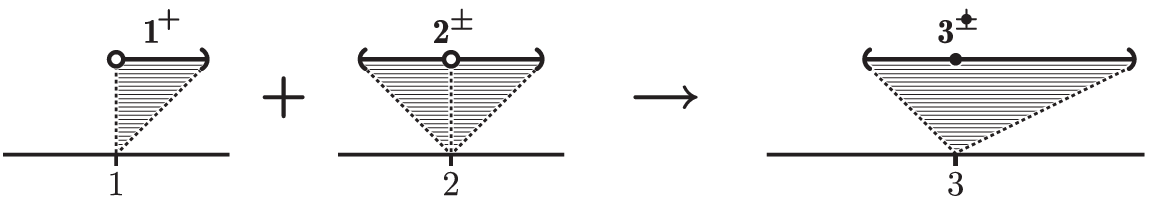}
\end{centering}
\caption{Addition:  $1^\NNp + 2^\NNpm \to 3^\NNpmd$}
\end{figure}

In our first example of near-number addition, we have $1^\NNp +\, 2^\NNpm \to 3^\NNpmd$.  Focusing on the top slices, which set the tone for our computation, we can quickly determine the sum, for the sum of two intervals is itself an interval, and the endpoints of the sum of two intervals are given by adding the respective endpoints of the summand intervals.  If we choose to compute in detail, we may assume the widths of the intervals to be any positive value $\varepsilon>0$ that we'd like, as they will shrink to zero as we move downward; $(1,1+\varepsilon)+(2-\varepsilon,2) = (3-\varepsilon,3+\varepsilon)$ and $(1,1+\varepsilon)+(2,2+\varepsilon)=(3,3+2\varepsilon)$, so our result is the union of these two intervals, $(3-\varepsilon,3+2\varepsilon)$.  Of course, we needn't go into this level of detail if our only concern is the basic shape of the result, but the ease of interval addition here makes it a worthwhile exercise.

\begin{figure}[!h]
\begin{centering}
\epsfbox{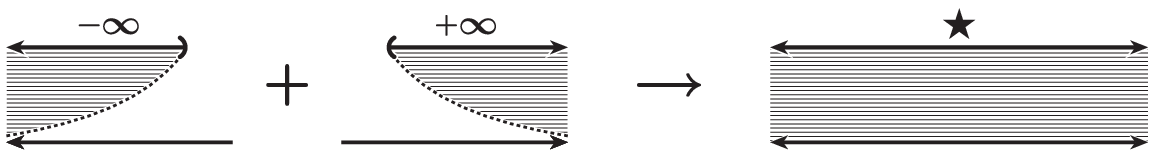}
\end{centering}
\caption{Addition:  $-\infty + \infty \to \NNstar$}
\end{figure}

In our second example of addition, $-\infty + \infty \to \NNstar$, we can clearly \emph{see} the indeterminacy of the sum.  With explicit slices of $-\infty$ and $+\infty$ in hand, we can choose any value in the left slice and any value in the right, obtaining some real number $A$ as their sum.  Moving the left-hand value abitrarily far leftward, we see that the result contains all values less than $A$; similarly, we can move the right-hand value arbitrarily far rightward, finding that the result contains all values larger than $A$.  Thus, the resulting slice is the entire real line at each level, yielding an indeterminate result for the sum.

\begin{figure}[!h]
\begin{centering}
\epsfbox{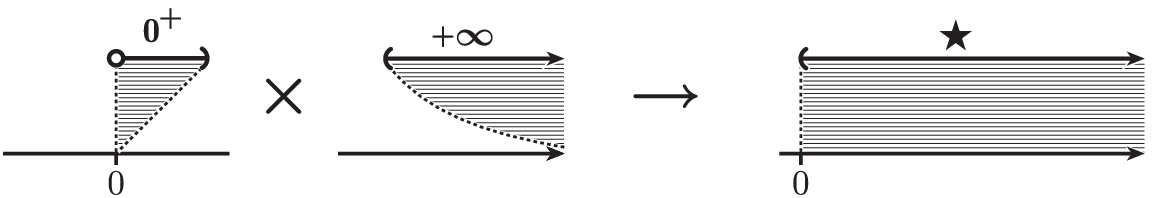}
\end{centering}
\caption{Multiplication:  $0^\NNp \times \infty \to \NNstar$}
\end{figure}

In our final example, we see the product $0^\NNp\times\infty\to\NNstar$, whose indeterminacy is again easily determined from its slices.  Taking slices of the factors, we can choose any value in the left-hand slice and any value in the right, obtaining some real number $A>0$ as their product.  Moving the left-hand value arbitrarily close to $0$, we see that the result contains all values in $(0,A]$; similarly, moving the right-hand value arbitrarily far rightward, we see that the result contains all values in $[A,\infty)$.  Thus, the resulting slice is the interval $(0,\infty)$, again yielding an indeterminate result. 

Note that in both of the previous cases of indeterminate near-number arithmetic, we have been able to see concrete diagrams demonstrating the indeterminacy of the results---thus we have more than memorization and intuitive explanation as rationale for such indeterminacy.  This geometric approach to near-number arithmetic provides the student a means to understand both usual arithmetic and indeterminate forms in a parallel fashion, as well as a means of rejustifying or re\"exploring any forgotten or troublesome case of the rules for arithmetic.  As such, we are provided a scaffold to support the arithmetic of near-numbers as their basic arithmetic rules are learned and internalized.

\subsection{Near-numbers and the limit}

We now consider the relationship between the expression $\lim_{x\to a}f(x)=L$ and the near-number system.  We say that this limit \emph{converges} just when $f(a^\NNpm) \to L^\NNpmd$ for some real number $L$; in this case, we call $L$ the \emph{limit}.  In the case that no such real number $L$ exists, we say that the limit \emph{diverges} (this includes the case that $f(a^\NNpm)\to\pm\infty$).  Should we simply know that $f(a^\NNpm)\to\NNstar$, we have an \emph{indeterminate} result, and we cannot \emph{a priori} conclude whether the limit converges or diverges---for even the finite near numbers $\alpha$ have $\alpha\to\NNstar$.  We must, in this case, either re-compute the limit via some other means (e.g., after algebraic manipulation) or, if we claim that it diverges, directly address exactly what $f(a^\NNpm)$ \emph{is}.  As an example of this last situation, we return to the function graphed Figure~\ref{Stack} and consider $f(2^\NNpm)\to 2^\NNm \cup 3^\NNm$.  While the only valid statement we can make in terms of the \emph{basic} near-numbers is ``$f(2^\NNpm)\to\NNstar$'', we can see from the diagram that this limit diverges, for the result fits no finite near-number (in each slice it retains values differing by more than $1$, a property that distinguishes it from any finite near-number).

Of course, the above discussion directly applies to all one-sided and infinite limits as well, simply replacing $x\to a^\NNpm$ with the appropriate $x\to\alpha$ in the limit expression; the near-number system is sufficiently robust to handle all limit cases seamlessly.  Certainly, there is a loss of information in converting the near-number result into a real number, and this is precisely how we should think of the classical limit operation: it \emph{crystallizes} a finite near-number into the real number at its base.  Such conversion into real numbers must only be performed at the very end of any computation involving limits, for as we've seen, the real numbers do not suffice to express intermediate results.

A near-number expression's \emph{structure} can greatly affect the result
when it is simplified, which accurately reflects our process of limit
computation.  
For instance, directly substituting for $x$ in the expression below, we obtain an 
indeterminate result:

\centerline{$\ds\left.\frac{x^2-x}{x^2 + x}\right|_{x\to\infty} \to 
\frac{\infty^2 - \infty}{\infty^2 + \infty} \to 
\frac{\infty - \infty}{\infty+\infty} \to 
\frac{\NNstar}{\infty} \to 
\NNstar$}

\vspace{6pt}

\noindent This result is not in conflict with what we know about $\lim_{x\to\infty}\frac{x^2-x}{x^2+x}$, because we've merely stated that the result $\to\NNstar$ (and, after all, every near-number $\to\NNstar$).  
Of course, if we factor and cancel (for $x\ne0$), we obtain the usual result:

%\centerline{$\ds\frac{x^2-x}{x^2 + x} = \frac{x^2(1 - \frac1x)}{x^2(1 + \frac1x)} = \frac{1 - \frac1x}{1 + \frac1x}$ (when $x\ne0$), and}

\vspace{6pt}

\centerline{$\ds\left.\frac{1 - \frac1x}{1 + \frac1x}\right|_{x\to\infty} \to
\frac{1 - \frac1{\infty}}{1 + \frac1{\infty}} \to
\frac{1 - 0^\NNp}{1 + 0^\NNp} \to
\frac{1^\NNm}{1^\NNp} \to
1^\NNm$}

\vspace{6pt}

\noindent%
The directionality of the $\to$ relation is the key component that allows us to obtain different results without contradiction; both results are valid, but the second is more precise.  The near-number system allows us to express our work in both successful and unsuccessful attempts at such a computation.  This is important, for even in a failed computation, writing out our work clearly provides us with
something concrete to examine in our search for a cue as to what method might
resolve the indeterminacy.

\def\sos{\mathrm{sos}}

As another example, we take the classical limit $\lim_{x\to0}\frac{\sin x}x=1$, which we could write as $\sos(0^\NNpm)\to1^-$, where $\sos(x) = \frac{\sin x}x$.  Nothing in our mental approach to such limits changes under the near-number system, but we are again prepared to show every step.  For instance:

\vspace{6pt}

\centerline{
$\ds\left.\frac{\sin(4e^x)}{e^x}\right|_{x\to-\infty} \to
\frac{\sin(4e^{-\infty})}{e^{-\infty}} \to
\frac{\sin(4\cdot 0^\NNp)}{0^\NNp} \to
\frac{\sin(0^\NNp)}{0^\NNp} \to
\frac{0^\NNp}{0^\NNp} \to
\NNstar$
}

\vspace{6pt}

\noindent
Having been able to write out our work, we are directly led to the proper method even if it was not apparent at the start.  Backtracking one step at a time from the $\NNstar$ to $\frac{\sin(0^\NNp)}{0^\NNp}$, we are cued to look for the $\sos$ function in our expression:

\vspace{6pt}

\centerline{
$\ds\frac{\sin(4e^x)}{e^x}
= \frac44\cdot\frac{\sin(4e^x)}{e^x}
= 4\cdot\frac{\sin(4e^x)}{4e^x}
= 4\cdot\sos(4e^x)$,}

\vspace{2pt}
\centerline{and}
\vspace{2pt}

\centerline{$4\cdot\sos(4e^{-\infty}) \to 
4\cdot\sos(4\cdot0^\NNp) \to 
4\cdot\sos(0^\NNp) \to 
4\cdot1^\NNm \to
4^\NNm$}

\vspace{6pt}

The near-number system allows us to fully express our work in limit computations, opening the process to full inspection by both student and teacher.  It provides
a precise and unambiguous language for written expression, removing many
of the obstacles to clear understanding that occur when the discussion is forced into the 
language of real numbers and equality.  We thus see a few of the linguistic benefits of employing the system of near-numbers in the computation of limits.

%\subsection{Formalizing arrow statements:  quantifying near-numbers}
\subsection{Quantifying near-numbers and the arrow relation}
\label{Quantify}

Before we're prepared to write a formal proof of a near-number statement, we must quantify our near-number diagrams.  After this simple process, we will easily be able to use our geometric understanding of the arrow relation to write symbolic statements ready for proof.  Here, we take the static planar-set view of near-numbers, as we desire a view of all slices of a near-number simultaneously.

\begin{figure}[!h]
\epsfbox{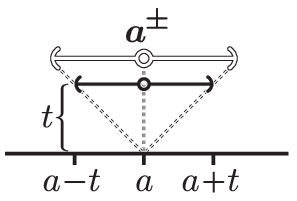}
\hspace{0.5in}
\epsfbox{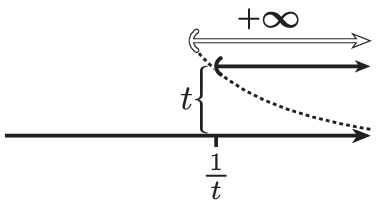}
\caption{Quantifying the near-numbers}
\label{QuantifyFig}
\end{figure}

Labeling the horizontal slices of each near-number $\alpha$ according to their height $t > 0$, we consider each slice to give a subset of $\R$.
In this context, the value $t$ measures proximity to $\alpha$, for as
$t$ decreases, these slices tighten in on the meaning of the near-number 
$\alpha$, as in Figure~\ref{QuantifyFig}.

\noindent%
With this in mind, we 
call a real number \emph{$t$-near} $\alpha$ if it lies within the 
$t$-slice of $\alpha$.  We can now read our near-number diagrams quantitatively
to obtain the conditions for being $t$-near each type of basic near-number, as in Figure~\ref{Intervals}.

\def\figbox#1{\mbox{\hfill\makebox[2in]{\hfill\epsfbox{#1}\hfill}\hspace{0.3in}}}

\begin{figure}[!h]
\begin{minipage}{0.9\textwidth}

\hspace{0.8in}\epsfbox{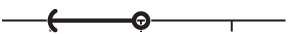}\hspace{0.2in}$a^\NNm:\;\;a-t < \square  < a$

\hspace{0.8in}\epsfbox{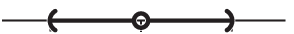}\hspace{0.2in}$a^\NNpm:\;\;a-t < \square  < a+t, \;\square\ne a$

\hspace{0.8in}\epsfbox{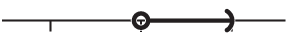}\hspace{0.2in}$a^\NNp:\;\;a < \square < a+t$

\vspace{8pt}

\hspace{0.8in}\epsfbox{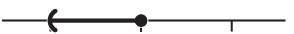}\hspace{0.2in}$a^\NNmd:\;\;a-t < \square \le a$

\hspace{0.8in}\epsfbox{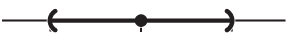}\hspace{0.2in}$a^\NNpmd:\;\;a-t < \square < a+t$

\hspace{0.8in}\epsfbox{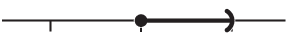}\hspace{0.2in}$a^\NNpd:\;\;a \le \square < a+t$

\vspace{-3.5pt}

\hspace{0.8in}\epsfbox{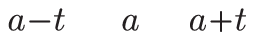}
\vspace{8pt}

\def\figbox#1{\hbox{\hfill\makebox[2.75in]{\hfill\epsfbox{#1}\hfill}\hspace{0.3in}}}

\epsfbox{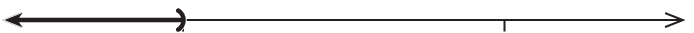}\hspace{0.2in}$-\infty:\;\; \square < -\frac1t$

\epsfbox{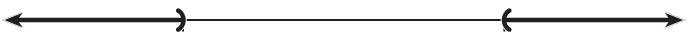}\hspace{0.2in}$\pm\infty:\;\; \square < -\frac1t$ or $\square > \frac1t$

\epsfbox{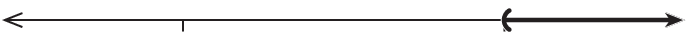}\hspace{0.2in}$+\infty:\;\; \square > \frac1t$

\vspace{-4pt}
\figbox{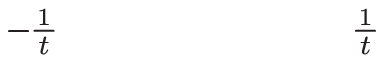}
\end{minipage}

\caption{Conditions for being $t$-near the basic near-numbers}
\label{Intervals}
\end{figure}

With our slices labeled, we are ready to translate our geometric definition of any basic near-number assertion $\alpha\to\beta$ into a formal logical definition as follows.  Recall that $\alpha\to\beta$ can be geometrically interpreted as ``$\alpha$ fits into $\beta$''---i.e., dropping the slices of $\alpha$ into the shape of $\beta$, the slices reach the bottom of $\beta$.  More specifically:

\noindent%
\begin{minipage}{0.75\textwidth}
\vbox{
\vspace{6pt}
\hangindent=0.5in
\hangafter=0
For each slice of $\beta$, 
some slice of $\alpha$ fits inside it.

\vspace{6pt}
Labeling the slices in $\alpha$ via $\delta>0$ and those in $\beta$ via $\varepsilon>0$:

\vspace{6pt}
\hangindent=0.5in%
\hangafter=0
For each $\varepsilon>0$, there is some $\delta>0$ so that 

\hangindent=0.5in%
\hangafter=0
the $\delta$-slice of $\alpha$ fits within the $\varepsilon$-slice of $\beta$.
\vspace{6pt}

Interpreting one slice fitting inside another as every value of the first lying in the second (an if-then statement), we obtain the general form of the definition:  

\vspace{6pt}
\fbox{\vbox{
\hangindent=0.5in%
\hangafter=0
For each $\varepsilon>0$, there is some $\delta>0$ so that 

\hangindent=0.5in%
\hangafter=0
if [$x$ is $\delta$-near $\alpha$] then [$x$ is $\varepsilon$-near $\beta$].
}}
\vspace{6pt}

}
\end{minipage}
\begin{minipage}{0.25\textwidth}
\vbox{
\hspace{0.25in}\epsfbox{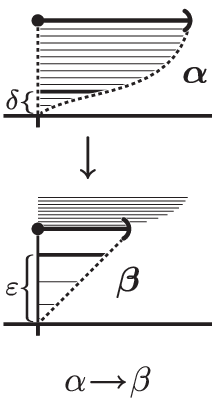}
\vspace{0.2in}
}
\end{minipage}

\noindent From this general template, all that remains to be done in order to formally define the statement ``$\alpha\to\beta$'' is to fill in the nearness conditions for the particular $\alpha$ and $\beta$ in question.  This system of interchangeable parts enables us to build all our definitions from a single foundation, thereby stressing the commonality of all cases.  In every case, the concept is the same:  the slices of $\beta$ stipulate position on the real line increasingly strictly;  $\alpha\to\beta$ means that $\alpha$ satisfies all of these conditions.

We note that, should it be desired, we can recover the ``standard'' notation for the infinite cases by labeling the slices of the infinite near-numbers via $N=\pm\frac1t$, as appropriate.  The trade-off in doing so is that a more standard and slightly simpler formulation for the infinite cases is obtained, at the cost of the loss of the completely parallel treatment of all near-numbers.  This affords a three-fold pedagogical decision on the treatment of the infinite near-numbers: to treat all near-numbers similarly, then to re-state the infinite cases; to label the slices for the infinite cases differently from the start; or to maintain parallelism throughout the presentation.

Via this system, we can begin to explore logical definitions and $\varepsilon$-$\delta$ proofs at a very basic level, working our way up to more complicated cases.  For example, the algebraic
elements of limit proofs can be practiced via proving statements such as ``If $x$ is $\frac1{10}$-near $1^\NNpmd$, then $2x+3$ is $\frac15$-near $5^\NNpmd$,'' while the logical elements can be practiced with minimal intrusion of algebra by proving statements such as ``$2^\NNp \to 2^\NNpm$'' and ``$2^\NNpd \not\to 2^\NNpm$''.

In the near-number system, each definition stems directly from our simple geometric definition of the arrow relation, providing an intuitive scaffold for the logical definition.  We thus eliminate the necessity for rote memorization as an \emph{a priori} method of instruction of limits---the statement is no longer ``here is the correct definition, and here is what is wrong with each of its permutations,'' but rather, ``here is the definition, which clearly and naturally captures our intent.''

\subsection{Quantifying with functions and arithmetic}

If we directly apply our basic definition to the statement $f(\alpha)\to\beta$, we obtain:

\vspace{6pt}
\hangindent=0.5in
\hangafter=0
For each $\varepsilon>0$, there is some $\delta>0$ so that 

\hangindent=0.5in
\hangafter=0
the $\delta$-slice of $f(\alpha)$ fits within the $\varepsilon$-slice of $\beta$.
\vspace{6pt}

\noindent%
Working with the graph of $f$, we can use the animated perspective on near-numbers 
to analyze this statement, placing $\alpha$ on the horizontal axis
and $\beta$ on the vertical axis.  Each $\varepsilon$-slice of $\beta$ paints a 
target, and for each such target there must be some 
$\delta$-slice of $f(\alpha)$ that fits into it.  

\begin{figure}[!h]
\epsfbox{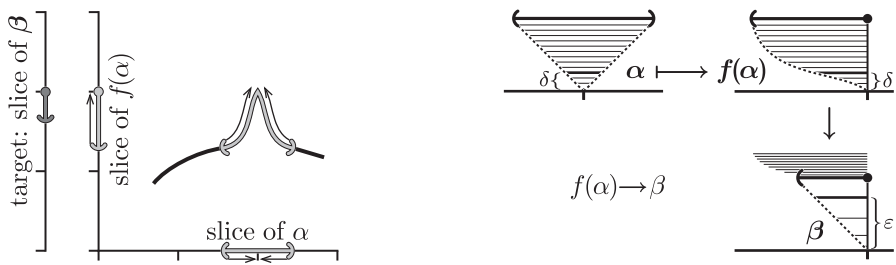}
\caption{$f(\alpha)\to\beta$:  animated and static-set perspectives}
\end{figure}

\noindent%
Symbolically, we see that the only difference between 
the definition of $\alpha\to\beta$ and $f(\alpha)\to\beta$ is that we map values
in $\alpha$ via $f$ before checking them to be within $\beta$ (via both our graphical analysis and the definition of $f(\alpha)$).  Thus, our symbolic definition is:

%\noindent%
%But the $\delta$-slice of $f(\alpha)$ is obtained by applying $f$ to the
%$\delta$-slice of $\alpha$, so this is:
%
%\vspace{6pt}
%\hangindent=0.5in
%\hangafter=0
%For each $\varepsilon>0$, there is some $\delta>0$ so that
%
%\hangindent=0.5in
%\hangafter=0
%$f[$the $\delta$-slice of $\alpha]$ fits within the $\varepsilon$-slice of $\beta$.
%\vspace{6pt}
%
%\noindent%
%To check this containment, we take a value $x$ in the $\delta$-slice 
%of $\alpha$, hit it with $f$, and check whether the resulting value $f(x)$ 
%lies in the $\varepsilon$-slice of $\beta$:

\noindent%
\begin{minipage}{0.75\textwidth}
\vbox{
\vspace{6pt}
\fbox{\vbox{
\hangindent=0.5in%
\hangafter=0
For each $\varepsilon>0$, there is some $\delta>0$ so that 

\hangindent=0.5in
\hangafter=0
if [$x$ is $\delta$-near $\alpha$] then [$f(x)$ is $\varepsilon$-near $\beta$].
}}
\vspace{6pt}

}
\end{minipage}

\noindent%
Filling in the relevant conditions for $\alpha$ and $\beta$, we can
immediately obtain the definition in any standard case of the limit.  
For example, $f(1^\NNpm)\to+\infty$ becomes:

\vspace{6pt}
\hangindent=0.5in
\hangafter=0
For each $\varepsilon>0$, there is some $\delta>0$ so that 

\hangindent=0.5in
\hangafter=0
if [$1-\delta < x < 1+\delta$ and $x\ne1$] then [$f(x) > \frac1\varepsilon$],
\vspace{6pt}

\noindent%
and $f(-\infty)\to 3^\NNpd$ becomes:

\vspace{6pt}
\hangindent=0.5in
\hangafter=0
For each $\varepsilon>0$, there is some $\delta>0$ so that 

\hangindent=0.5in
\hangafter=0
if [$x < -\frac1\delta$] then [$3\le f(x) < 3+\delta$].
\vspace{6pt}

The near-number system is flexible enough to seamlessly handle multivariable functions and arithmetic operations, as well.  For example, recalling that each point in the slice of the sum of two near-numbers is obtained by adding points of the corresponding slices of the summands, we can work from basic definitions to quantify $\alpha_1+\alpha_2\to\beta$ as:

\vspace{6pt}
\hangindent=0.5in
\hangafter=0
For each $\varepsilon>0$, there is some $\delta > 0$ so that 

\hangindent=0.5in
\hangafter=0
if [$x_1$ is $\delta$-near $\alpha_1$] and [$x_2$ is $\delta$-near $\alpha_2$] 

\hangindent=0.5in
\hangafter=0
then [$x_1+x_2$ is $\varepsilon$-near $\beta$].
\vspace{6pt}

Via the near-number system, we are able to discuss, quantify, and come to understand our objects of study via a structured approach:  first, the near-numbers themselves; second, the arrow relation for near-numbers; third, the arrow relation as applied to functions.  Thus, we can build our knowledge of the limit in
steps, rather than swallowing an atomic definition whole---we obtain the 
same logical conclusion from 
a progression of easily-understood geometric situations that naturally give
rise to the relevant formal definitions.

\section{Conclusion}

The near-number system provides the student a graded path via which to ascend the mountain of limits in the calculus curriculum---intuitively, computationally, and formally---so that during the learning process, smaller, clearer, and more intuitively motivated steps can be made.  Near-numbers provide a proper set
of \emph{objects} via which the relevant concepts can be discussed,
repairing the linguistic flaws inherent in the real-number treatment of limits
and avoiding the ambiguities and apparent contradictions that often cloud the student's view of limits.  At the logical peak of the mountain, near-numbers provide a clear means of reasoning out the proper definition of limit, and moreover, they do so seamlessly for all cases of one-sided, two-sided, and infinite limits.  Thus, not only does the limit lose its burden of atomic definitions in need of memorization, but the inherent similarity of the definitions in all limit cases is made evident.  Finally, the purpose of the limit (as a real number) is very apparent to the student:  it is to crystallize and make explicit certain notions objects that by their very nature are at first quite elusive.

\appendix
\section*{Appendix:  Mathematical formalization of the near-number system}

The near-number can be treated formally as a mathematical object in standard set theory, as will be described below.  The purpose of this section is not to suggest that a fully formal treatment of near-numbers need be introduced into the curriculum, any more than it would be suggested to formally define the real numbers via Dedekind cuts of rationals in an introductory calculus course.

Formally, a near-number is a function $\alpha:(0,1)\to \mathcal{P}(\overline{\R})$ for which $s \le t \Rightarrow \alpha(s) \subseteq \alpha(t)$.  Here, $\overline{\R} = \R \cup \lbrace \frownie \rbrace$, where the additional value $\frownie$ indicates an invalid computation (such as division by zero).  As special cases, we have the finite near-numbers such as $a^\NNp:t\mapsto(a,a+t)$ and the infinite near-numbers such as $+\infty:t\mapsto(\frac1t,\infty)$.  Within this context, each real number $a$ corresponds to the near-number $a^\NNdot:t\mapsto\lbrace a \rbrace$.  The indeterminate near-number $\NNstar$'s function as a catch-all is formalized by setting $\NNstar:t\mapsto\R\cup\lbrace\frownie\rbrace$.  To be more specific, we could call $\NNstar$ the \emph{exceptional near-number}, distinct from the merely indeterminate near-number $\Omega:t\mapsto\R$; however, in the interests of simplicity, we do not make this distinction in our presentation.

A real-valued function can be made to act on a near-number, in much the same way that it acts on a set of real numbers.  For any real-valued function $f$, we define a new function $\overline{f}$ by $\overline{f}(x) = f(x)$ if $x\in\mathrm{domain}(f)$ and $\overline{f}(x)=\frownie$ otherwise (which allows our function to act on any set, flagging any invalid arguments via $\frownie$).  Real-valued functions of two (or more) real variables, such as the arithmetical operations, can act on near-numbers similarly via $f(\alpha,\beta):t\mapsto\overline{f}[\alpha(t)\times\beta(t)]$  (while the symmetry in the parameters might catch the eye as problematic, it is not).  Finally, the near-number relation $\alpha\to\beta$ is defined as follows for near numbers $\alpha$, $\beta$: $\alpha\to\beta$ if $\forall\,t\in(0,1)\;\exists\,s\in(0,1)$ such that 
$\alpha(s)\subseteq\beta(t)$ (i.e., $\alpha$ ``fits inside'' $\beta$, as discussed in Section~\ref{Arrow}).

In contrast to some other alternative approaches to the limit, the notion of proper infinitesimals is abandoned in favor of the use of sequences of sets of real numbers to represent the key objects in the concept of limit, reflecting the standard formal nature of such objects.  The careful reader may note that the menagerie of near-numbers includes much more than just the basic ones first 
introduced.  For example, for any $\alpha,\beta$, we have their intersection $\alpha\cap\beta:t\mapsto\alpha(t)\cap\beta(t)$ and union $\alpha\cup\beta:t\mapsto\alpha(t)\cup\beta(t)$.  In addition, any set $S\subseteq\R$ gives a constant near-number $t\mapsto S$; e.g., $\sin(+\infty)\to[-1,1]$, which allows us to write:

\vspace{6pt}
\hfill $\ds\frac1{\infty}\sin(\infty) \to 0^\NNp\times[-1,1] \to 0^\NNpmd$ \hfill
\vspace{6pt}

\noindent While not a substitute for the full ``Squeeze Theorem'', near-numbers do give us the tools to correctly write down what we think when we see the limit $\lim_{x\to\infty} \frac1x\sin x$.

We remark that much greater symmetry is obtained in the near-number system 
if we compactify the real line to a circle as usual by adding a point $\infty$,
obtaining six near-numbers $\infty^\NNp$, $\infty^\NNm$, $\infty^\NNpm$, $\infty^\NNpd$, $\infty^\NNmd$, and $\infty^\NNpmd$ instead of our three special infinite near-numbers---in this case, all basic near-numbers other than $\NNstar$ have the
same ``shapes''.  Under the convention that the reciprocal of $0$ is $\infty$ (treating this operation geometrically as reflection, rather than as
multiplicative inversion),
we can then dispense with all special cases for reciprocation.  The benefit of
the near-number system in this context is that the conflict between the desire
for a single point at infinity and the need for distinct limits $\lim_{x\to+\infty}$ and $\lim_{x\to-\infty}$ is resolved, for $\infty^\NNp$ and $\infty^\NNm$ smoothly treat the two sides of the single value $\infty$ as needed (although the usual issues of operation with the value $\infty$ itself remain at issue).

\end{document}